\documentclass [12 pt]{article}
\usepackage[english]{babel}
\usepackage[T1]{fontenc}
\usepackage[latin1,utf8]{inputenc}
\usepackage[top=72.27pt,bottom=72.27pt, left =72.27pt, right=72.27pt]{geometry}
\usepackage{lmodern}
\usepackage{amssymb,amsmath}
\usepackage{type1cm}
\usepackage{latexsym}
\usepackage{graphicx}
\usepackage{color}
\usepackage{float, caption}
\usepackage{titlesec}
\usepackage{tikz}
\usepackage{url}
\usepackage{multirow}

\usepackage{pdfpages}
\pdfoutput=1

\titleformat{\section}{\large\bfseries}{\thesection.\!\! \!}{1em}{}

\begin{document}
{

\title{\textbf{Integrating Storage to Power System Management }}
\author{ Luckny Zéphyr, C. Lindsay Anderson
\\
\small{Biological and Environmental Engineering Department, Cornell University}
}
\maketitle

\begin{abstract}
Wind integration in power grids is very difficult, essentially because of the uncertain nature of wind speed. Forecasting errors on output from wind turbines may have costly consequences. For instance, power might be bought at highest price to meet the load. On the other hand, in case of surplus, power may be wasted. Energy storage facility may provide some recourse against the uncertainty on wind generation. 
 Because of the sequential nature of power scheduling problems, stochastic dynamic programming is often used as solution method. However, this scheme is limited to very small networks by the so-called curse of dimensionality. To face such limitations, several approximate approaches have been proposed. We analyze the management of a network composed of conventional power units as well as wind turbines through approximate dynamic programming. We consider a general power network model with ramping constraints on the conventional generators.  We use generalized linear programming techniques to linearize the problems. We test the algorithm on several networks of different sizes and report results about the computation time. We also carry out comparisons with classical dynamic programming on a small network. The results show the algorithm seems to offer a fair trade-off between solution time and accuracy.
\end{abstract}

\section{Introduction}
Cost-effective management of power units is a very challenging task. Generating units have to be committed such that demand for electricity is met. This is more difficult to carry out  because of significant load variability, among others. 

Commitment of traditional units such as nuclear and coal-fire involve cost, mostly because of fuel consumption, and is source of environment concerns. Driven by both increasing environmental awareness and technological advances, since the last decade, wind-based electricity generation has been widely promoted \cite{zhou2014managing}. For instance, the 2001/77/EC European Commission Directive had set at 22\% the  renewable integration target for Europe by 2010 \cite{moura2010role}. In the United States, in 2006, the total wind installed capacity grew approximately from 9 000 MW to 11 600 MW \cite{smith2007utility}. As a result of  political emphasis, in 2006, on the need to increase the United States energy efficiency and to diversify the energy portfolio, a collaborative initiative was created to explore the requirements for a 20\% wind integration scenario by 2030 \cite{lindenberg200920}.

However, unlike conventional power source, output from wind turbines is uncontrollable. Due to the random nature of such power, these latter cannot serve both supply and reliability purposes (in case of outage in the network). As a consequence, wind generators cannot replace conventional ones to meet peak load \cite{moura2010role}.  The intermittency of wind generation may create an imbalance between the supply and demand for power. In case of excess generation, power may be curtailed due to transmission congestion, and in case of peak load, to meet the demand, power would probably be bough at expensive price because of inadequate anticipation (forecasting error).

 Energy storage devices may serve as recourse to circumvent the uncertainty of wind power.  Such devices may be used to store excess generation, or for arbitrage profits. Indeed power might be bough at a lower price during off-peak hours to be stored and sold at a higher price during peak load hours \cite{mokrian2006stochastic}. For an in depth description of existing and in development storage technology, see \cite{ibrahim2008energy, hadjipaschalis2009overview}. Benefits as well as market-related questions are discussed in \cite{eyer2004energy}. 

A significant stream of research has been focusing on harvesting wind power in the presence of storage. The benefits of coupling intermittent renewable energy with storage are discussed in \cite{barton2004energy, mcdowall2006integrating, barton2006probabilistic, black2006value}. Prior studies have also analyzed the coupling of wind generation with storage via stochastic programming (e.g. \cite{garcia2008stochastic}, \cite{abbey2009stochastic}, \cite{meibom2011stochastic}). Unlike other research where a fixed wind generation curve is used, \cite{succar2012optimization} analyze the optimal wind turbine rating and the storage configuration of a wind farm coupled to compressed air energy storage.  \cite{scott2012approximate} compare approximate dynamic programming schemes for energy storage management based on instrumental Variables and projected Bellman errors for a model ignoring transmission lines in the presence of a single storage device. In \cite{zhou2014managing}, the valuation of  a single storage unit provided by a stochastic approach is compared with several valuation heuristics for a model with limited transmission capacity. Storage of energy for arbitrage opportunities is studied by \cite{mokrian2006stochastic} via stochastic programming and dynamic programming.

We examine the management of a power network composed of conventional units, wind turbines, as well as storage devices. In particular, we discuss an approximate stochastic method based on stochastic dual dynamic programming for the management of the storage units. Section \ref{sec:description} deals with the operation of a power network equipped with energy storage facilities. Section \ref{sec:SDP} provides a formulation of the problem under the framework of stochastic dynamic programming and discusses the associated sources of complexity.  Section \ref{GLP} discusses a linear approximation of the problem. In Section \ref{sec:ASDP}, we review some general approximate dynamic programming techniques widely used in power system management, while Section \ref{sec:SDDP} presents an approximate formulation for the problem. We illustrates the algorithm's performance over a small network and carry out comparisons with dynamic programming in Section \ref{sec: sddp-dp-illus}. Results of  additional experiments on larger networks are reported in Section \ref{sec:more-num}.

\section{Operation and control of a  power system with storage} \label{sec:description}

Optimization models are often used to determine optimal \textit{generation schedule} for the conventional generators over a finite planning horizon, $T$, in order to meet the demand in each time step $t$ at minimal cost. This is usually achieved by formulating the problem as a two-stage optimization problem, where in the first stage decisions are made about on/off status of the conventional generators based on anticipation (or forecasting) on demand and wind generator output. The second stage decisions are usually concerned with generation (\textit{economic dispatch} of the committed generators), as well as recourse actions (reserve requirements) to offset forecasting errors and for reliability purposes in case of unforeseen events such as generator failure, transmission line disruption. The time step is usually an hour or a fifteen-minute interval, and the planning span is usually 24, 48 or 168 hours (a week). 

Based on load and wind generation forecasts, the operator makes decision about the commitment of the conventional generators in each time period of the horizon. We assume that generation decisions (for  conventional generators) are made after observations of the random variables (wind generator output).

 Here, we are dealing with the economic dispatch problem alone, assuming that the conventional generators are all committed. For instance, \cite{asamov2015regularized} first use a simulator to determine the on/off status of the traditional generators in each point of time, then solve the dispatch problems.

A power system may be represented by a directed graph $\mathcal{G}=(N,L)$, where each node $n\in N$ represents a bus, and each link a transmission line.  In each node  may be located a set of conventional generators $G_n$, a wind farm, and a storage device. We distinguish the different components of the graph via the following set of indices: (i) $G$ is the set of conventional generators, (ii)  $M \subseteq N$ is the set of wind farms, and (iii) $S$  denotes the set of storage devices.

The generators are mechanical devices, they have to operate within finite generating limits.  While in motion, a minimal output may be required for a generator to be in steady state. Similarly, a threshold is imposed on the maximal output to avoid mechanical damages. For any generator $g \in G, \underline{p}_g$ and $\bar{p}_g$ (in $MW$) denote these bounds, respectively. In addition, there usually are minimum and maximum allowable changes in the generations between two consecutive  periods. For any generator $g \in G$, the following defines ramp down and ramp up limits:
\[-\underline{\lambda}_g \leq p_{gt}-p_{g,t-1}\leq \bar{\lambda}_g.\]

In addition to power generated in a node, power can flow between two nodes trough transmission lines defined as the set of undirected pairs $L\subseteq N \times N$  (assuming power can flow in either direction). For any node $n\in N$, define $O_n=\{(n,j)\in L\}$ to be the set of transmission lines that leave node $n$. Similarly, let $I_n=\{(j,n)\in L\}$ be the set of transmission lines that enter node $n$.  

For any node $n$, $\hat{D}_{nt}$  denotes a ``reliable'' forecast (in $MW$) of the load (demand),  in period  $t$. This demand is supplied by the generation in the node, by power transported to the node through the transmission lines, or by discharging the storage, if any, or by some combination of power from the three.  Thus, $p_{gt}$ denotes power generation, in $MW$, from generator $g$ in period $t$ at the cost $CP_{gt} (p_{gt})$. 

The output from the wind turbines is random. Indeed, generation from such units is conditioned on wind speed, which depends on uncontrollable meteorological conditions. Therefore, $\tilde{W}_t=(\tilde{W}_{1t},\dots, \tilde{W}_{|M|t})$ denotes the random vector of outputs from the wind turbines in period $t$; $w_{mt}$, in $MW$, is a particular realization of the stochastic process $\{\tilde{W}_{mt}\}$.

 We note by $e_{lt}$ the power, in $MW$, flowing through transmission line $l$ in period $t$.  The transmission lines have limited capacity, thus upper bounds are imposed on the power flowing through to prevent disruption. $\bar{e}_l, l\in L$, denote such bounds. We assume bidirectional lines, i.e., power can flow in either direction of the line. Therefore, for any transmission line $l=(n,n')\in L$, in period $t$, if $e_{lt}$ is positive, power is transmitted from node $n$ to node $n'$. A negative value indicates the opposite. The power flowing through a transmission line is proportional to the difference between the phase angles, in radians, of the two end buses, i.e., $e_{lt}=B_{l}(\theta_{n't}-\theta_{nt})$, where $B_l$ is the susceptance of line $l$ \cite{papavasiliou2013multiarea}. To avoid that the system is over-determined, i.e. with more equations than unknown, the voltage angle at the reference (slack) bus is usually set to zero.

In each period $t$, in addition to generation decisions, the operator of the network also makes decisions on the use of the storage. In each node where a storage device is located, in each period, the device may be charged or discharged to compensate for imbalance between the supply and demand. Such imbalance is more likely to occur in peak demand hours, i.e., periods where electricity consumption is highest. For instance, demand for power tends to be lower at night when people are asleep. During the summer, peak demand tends to occur in the afternoon because of higher activity rate and use of air conditioning by buildings. In the winter, demand is often higher in the morning and in the evening due to greater need for heating. Operators usually prepare for peak demand by committing extra power plants that may be called upon quickly in periods of higher demand or failures (reserve).

For any storage device $j\in S$, define $s_{jt}$ to be its level of energy ($MWh$) in the beginning of period $t$ (or the end of period $t-1$). This energy is converted into power (MW) via the simple equation $\text{Energy}=\text{Power}\times \text{Time}$. We note by $\Delta_{jt}$ the variation in the level of charge from the beginning through the end of period $t$. A negative value of $\Delta_{jt}$ indicates discharging, and a positive value the opposite. There may be a cost associated with varying stored energy level. $CS_{jt}(\Delta_{jt})$ is the cost associated with such variation, if any. Such cost is proportional to the amount of energy stored or discharged \cite{schoenung2011energy}.

The storage units have limited capacity, and in order to last, they cannot be completely depleted. Consequently, the level of charge can only varied progressively over time. We therefore note by $\underline{s}_j$ and $\bar{s}_j$ lower and upper bounds on the level of storage of device $j$. Assume that the discharge and charge maximum capacity is the same, in period $t$, the variation in the level of device $j$ then obeys
\[-\underline{\Delta}_j\leq \Delta_{jt}\leq \bar{\Delta}_j.\]

Lastly, each storage facility has input and output efficiency; energy is lost both at charging and discharging. We then define $c_j$ (resp. $d_j$) to be the efficiency coefficient of charging (resp. discharging) of storage device $j\in S$, where $0<c_j\leq 1$, and $0<d_j\leq 1$.

Imbalance between the load and supply may result because of wind power forecasting errors.  Since the storage units have limited capacity, they may not be able to absorb the total production surplus (charging) or to deliver the difference (discharging) in case of power shortage. For any node $n \in N$, $\gamma_{nt} \in \bold{R}$ denotes the power absorbed or delivered by the storage unit. In case of imbalance, we assume that  power excess (shortage) is absorbed (delivered) at a high rate. For each node $n\in N, \kappa^+_{nt} \in \bold{R^+}$ ($\kappa^-_{nt} \in \bold{R^+}$) denote such excess (shortage); $\bold{R}$ and $\bold{R^+}$ are the set of real and non-negative real numbers, respectively.

 Let us also explicitly distinguish charging and discharging via the variables $\Delta^+_{jt}\in \bold{R^+}$, and $\Delta^-_{jt}\in \bold{R^+}, j \in S$, respectively. Thus, in each period $t$, we have the power balance equations
\[\sum^{N_n}_{k=1}p_{kt}+\sum_{l\in I_n}e_{lt}+w_{nt}-\delta_n\gamma_{nt}-\kappa^+_{nt}+\kappa^-_{nt}=\sum_{l\in O_n}e_{lt}+\hat{D}_{nt}, \forall n\in N,\]
with 
\begin{equation*} 
\delta_{n}=\left\{ \begin{array}{ll}
1 & \text{if  there is a storage facility at node n},\\
0 & \text{otherwise}.
\end {array} \right.
\end{equation*}
By convention, $w_{nt}=0, n\in N$, if there is no wind farm at node $n$. We also have
\[\Delta^+_{jt}=\max\{0,c_j\gamma_{jt}\},  j \in S, \] and
\[\Delta^-_{jt} = \max\left\{0, - \frac{\gamma_{jt}}{d_j}\right\},  j \in S. \]
The storage level then evolves according to
\[s_{j,t+1}=\alpha_j s_{jt}+\Delta_{jt}^+ - \Delta_{jt}^-, \forall j\in S,\]
where $\alpha_j \in (0,1]$ is the storage efficiency of device $j$. As defined, it is obvious that both charging and discharging will not occur at the same time.

We aim to find the \textit{policy} $(P^*,\Delta^*)=\left[(P_{1}^*,\Delta_{1}^*),\cdots, (P_{T}^*,\Delta_{T}^*)\right]$ that minimizes the expected operating cost over the entire planning horizon, where $P_{t}^*=(p_{1t}^*,\cdots, p_{|G|t}^* )$, and $\Delta_{t}^*=(\Delta_{1t}^*,\cdots, \Delta_{|S|t}^* ), 1\leq t\leq T$; $\Delta_t$ refers to either charging or discharging.  $(P^*,\Delta^*)$ then solves

\begin{align}
&\underset{P_t,\Delta_t,\kappa_t,s_{t+1}}{\min{}}\left\{\mathbb{E}\left[\sum_{t=1}^T\left(\sum_{g=1}^{|G|}CP_{gt} (p_{gt})+\sum^{|S|}_{j=1}CS_{jt}(\Delta_{jt}) + M \sum ^N_{n=1}(\kappa^+_{nt}+\kappa^-_{nt}) \right)\right]\right\} \label{eq:ProbGenSt}\\
& \text{Subject to, for  $1\leq t\leq T:$}\notag \\
& \qquad  \qquad \quad \sum^{N_n}_{k=1}p_{kt}+\sum_{l\in I_n}e_{lt}+w_{nt}-\delta_n\gamma_{nt}-\kappa_{nt}^+ +\kappa_{nt}^-=\sum_{l\in O_n}e_{lt}+\hat{D}_{nt}, \quad \ \ \quad  n\in N \label{eq:ProbGenConstSt}\\
& \qquad  \qquad \quad  p_{gt}-p_{g,t-1}\leq \bar{\lambda}_g, \qquad \quad \qquad \ \qquad \quad \qquad \qquad \quad \qquad  \qquad \qquad  \quad g\in G \label{ineq:boundsptstart}\\
& \qquad  \qquad \quad p_{g,t-1} -  p_{gt} \leq \underline{\lambda}_g, \qquad \quad \qquad \ \qquad \quad \qquad \qquad \quad \qquad  \qquad \qquad  \quad g\in G\\
& \qquad  \qquad \quad \underline{p}_{g}\leq p_{gt} \leq \bar{p}_g, \quad \qquad \quad \quad \qquad \quad \quad \qquad  \qquad \ \qquad \quad \qquad \qquad \quad g \label{ineq:boundsptend} \in G\\
& \qquad  \qquad \quad  e_{lt}=B_{l}(\theta_{n't}-\theta_{nt}), \quad \quad \qquad \quad\quad \qquad \qquad \quad \qquad \quad \qquad \qquad \quad l=(n,n')\in L\\
& \qquad  \qquad \quad  -\bar{e}_l \leq e_{lt}\leq \bar{e}_l, \quad \quad \qquad \quad\quad\quad \quad \qquad \quad \qquad \quad \qquad \qquad \qquad \ l\in L \\
& \qquad  \qquad \quad  s_{j,t+1}=\alpha_j s_{jt}+\Delta^+_{jt} - \Delta^-_{jt},   \quad\quad \ \quad\quad \quad \qquad \quad \qquad \quad \qquad \qquad j \in S \\
& \qquad  \qquad \quad \underline{s}_j\leq s_{j,t+1}\leq\bar{s}_j, \quad \qquad \ \quad \qquad \qquad  \quad\quad\quad \qquad \ \quad \qquad \qquad  \quad j\in S \\
& \qquad  \qquad \quad 0\leq \Delta^+_{jt}\leq \bar{\Delta}_j, \quad \qquad \ \quad \qquad \qquad \quad\quad \qquad \quad \qquad \quad \quad \qquad  \quad \! j\in S \\
& \qquad  \qquad \quad 0\leq \Delta^-_{jt}\leq \bar{\Delta}_j, \quad \qquad \quad \qquad \quad \quad\quad \qquad \quad \qquad \quad \quad \quad \qquad  \quad \! j\in S\\
& \qquad  \qquad \quad \gamma_{jt}=\frac{\Delta^+_{jt}}{c_j} - d_j\Delta^-_{jt}, \ \quad \quad \qquad \qquad  \qquad \quad \qquad \quad \quad  \quad \quad  \quad \quad j\in S \label{eq:gamma}\\
& \qquad  \qquad \quad \kappa_{jt}^+ \geq 0,  \quad \qquad \quad \qquad \quad \quad  \quad \qquad \quad \qquad \quad \quad \quad \quad \! \qquad \quad \quad j\in S\\
& \qquad  \qquad \quad \kappa_{jt}^- \geq 0,  \quad \qquad \quad \quad \quad\quad \quad \quad\quad\quad \qquad \quad \quad  \quad \quad\ \qquad \quad \quad j\in S\\
& \qquad  \qquad \quad \Delta^+_{jt} \geq 0,  \quad\quad \qquad \quad \quad \quad \qquad \quad \qquad \quad \quad  \quad \quad \quad \! \qquad \quad \quad j\in S\\
& \qquad  \qquad \quad \Delta^-_{jt} \geq 0,  \quad \qquad \quad  \quad \qquad \quad \quad \quad\qquad \quad \quad  \quad \quad \quad \! \qquad \! \qquad j\in S  \label{eq:ProbGenEn}
\end{align}
$\mathbb{E}$ is the expectation operator, which is taken over $\tilde{W}_t$. For any set X, |X| denotes its cardinality. $M$ is a big number.

The multi-period problem (\ref{eq:ProbGenSt}--\ref{eq:ProbGenEn}) can theoretically be solved to optimality if functions $CP_{gt}$ is convex in $p_{gt}$,  $CS_{jt}$ being proportional to $\Delta_{jt}$. Indeed, the set 

\noindent $\Psi_t=\left\{(P_t, e_t,\tilde{W}_t,\theta_t,s_t,\Delta_t^+,\Delta_t^-,\gamma_t,\kappa_t^+,\kappa_t^-)|(\ref{eq:ProbGenConstSt}-\ref{eq:ProbGenEn})\right\}$ is a polyhedron. Therefore, by convexity of the cost functions, (\ref{eq:ProbGenSt}--\ref{eq:ProbGenEn}) is a convex problem. However, the problem may not be tractable numerically if we want a detailed representation for the underlying process of $\tilde{W}_t$.

\section{Representation under the framework of stochastic dynamic programming} \label{sec:SDP}
Problem (\ref{eq:ProbGenSt}--\ref{eq:ProbGenEn}) is a multi-period stochastic program. It also is a sequential decision problem. In fact, in each period $t$, the operator of the system observes the level of the stored energy as well as the generations in the previous period, and based on updated forecast for the wind turbines outputs (or wind speed) and the observation of the previous outputs, he/she decides about the generation of the conventional generators as well as the variation in the storage (charging or discharging) in order to meet the demand in each node. Thus, the tuple $(s_t,p_{t-1},w_{t-1})$ will be called the \textit{state of the system}, or state for short. 

Let us observe that whereas $p_{t-1}$ is known, we see from constraints (\ref{ineq:boundsptstart}-\ref{ineq:boundsptend}) that $p_{g,t-1}-\underline{\lambda}_g$ and $\underline{p}_g$ are two minorants for $p_{gt},  g \in G$. Similarly, $p_{g,t-1}+\bar{\lambda}_g$ and $\bar{p}_g$ are two majorants for $p_{gt}$. Consequently, we have
\[\max\{p_{g,t-1}-\underline{\lambda}_g,\underline{p}_g\}\leq p_{gt} \leq \min \{p_{g,t-1}+\bar{\lambda}_g, \bar{p}_g\}.\]

In principle, stochastic dynamic programming (SDP) is suited for problem (\ref{eq:ProbGenSt}--\ref{eq:ProbGenEn}). SDP sequentially decomposes (by period) the overall problem into smaller subproblems in a coordinated way, by seeking the best trade-off between the immediate and future use of the storage. Define $F_t (s_t, p_{t-1}, w_{t-1})$ to be the \textit{cost-to-go function} from the beginning of period $t$ to the end of the horizon. In addition, assume that the process $\{\tilde{W_t}\}$ is Markovian, i.e., $\mathbb{P}(\tilde{W}_t=w_t|\tilde{W}_{t-1}=w_{t-1},\ldots,\tilde{W}_0=w_0)=\mathbb{P}(\tilde{W}_t=w_t|\tilde{W}_{t-1}=w_{t-1})$. Therefore, for $t=T,T-1,\cdots,1$, a SDP recursion associated with problem (\ref{eq:ProbGenSt}--\ref{eq:ProbGenEn}) is given by
\begin{align}
F_t(s_{t}, p_{t-1}, w_{t-1})=&\underset{P_t,\Delta_t,\kappa_t,s_{t+1}}{\min{}}\left\{\sum_{g=1}^{|G|}CP_{gt} (p_{gt})+\sum^{|S|}_{j=1}CS_{jt}(\Delta_{jt}) + M \sum ^N_{n=1}(\kappa^+_{nt}+\kappa^-_{nt})+ \right. \notag\\
& \qquad \qquad \mathbb{E}_{\tilde{W}_{t}|w_{t-1}}\left[F_{t+1}(s_{t+1},p_t,\tilde{W}_{t})\right]\bigg\} \label{eq:SDPGenSt}\\
\text{S.t.} & \sum^{N_n}_{k=1}p_{kt}+\sum_{l\in I_n}e_{lt}+w_{nt}-\delta_n \gamma_{nt}- \kappa_{nt}^+ +\kappa_{nt}^-=\sum_{l\in O_n}e_{lt}+\hat{D}_{nt}, \ n\in N \label{eq:SDPGenConstSt}\\
&  p_{gt} \geq \max\{p_{g,t-1}-\underline{\lambda}_g,\underline{p}_g\},  \quad \quad \quad \quad\quad \quad \quad \quad\quad \quad\quad \! \quad \quad \quad \ g\in G \label{ineq:bIpt} \\
& p_{gt} \leq \min \{p_{g,t-1}+\bar{\lambda}_g, \bar{p}_g\},  \quad \quad \quad \quad\quad \quad \quad \quad\quad \quad\quad \! \quad \quad \quad \ g\in G \label{ineq:bSpt} \\
& e_{lt}=B_{l}(\theta_{n't}-\theta_{nt}), \quad \qquad \qquad \quad \qquad \quad \quad\quad \qquad \quad \quad l=(n,n')\in L\\
& -\bar{e}_l \leq e_{lt}\leq \bar{e}_l, \quad \qquad \quad \qquad  \qquad \qquad \quad \qquad \quad \qquad \qquad \! \ l\in L \\
&  s_{j,t+1}=\alpha_j s_{jt}+\Delta^+_{jt} - \Delta^-_{jt},     \quad \qquad \quad \ \quad \qquad \quad \qquad \quad \qquad \! j \in S \label{eq:stor}\\
&  \underline{s}_j\leq s_{j,t+1}\leq\bar{s}_j, \quad \qquad \ \quad \qquad \qquad  \quad\quad\quad \qquad \ \quad \qquad \! \quad j\in S \\
&  0\leq \Delta^+_{jt}\leq \bar{\Delta}_j, \!\quad \qquad \ \quad \qquad \qquad \quad\quad \qquad \quad \quad \quad \qquad  \quad \! j\in S \\
&  0\leq \Delta^-_{jt}\leq \bar{\Delta}_j, \!\quad \qquad \quad \qquad \quad \quad\quad \qquad \quad \quad \quad \quad \qquad  \quad \! j\in S\\
&  \gamma_{jt}=\frac{\Delta^+_{jt}}{c_j} - d_j\Delta^-_{jt}, \ \!\quad \quad \qquad \qquad  \qquad \quad \quad \quad  \quad \quad  \quad \quad j\in S \label{eq:gamma}\\
&  \kappa_{jt}^+ \geq 0,  \quad \qquad \quad \qquad \quad \quad \! \quad \qquad \quad \quad \quad \quad \quad \! \qquad \quad \quad j\in S\\
&  \kappa_{jt}^- \geq 0,  \quad \qquad \quad \quad \quad\quad \quad \quad\quad\quad \quad \quad  \quad \quad\ \qquad \quad \quad j\in S\\
& \Delta^+_{jt} \geq 0, \! \quad\quad \qquad \quad \quad \quad \qquad \quad \quad \quad  \quad \quad \quad \! \qquad \quad \quad j\in S\\
& \Delta^-_{jt} \geq 0,  \!\quad \qquad \quad  \quad \qquad \quad \quad \quad \quad \quad  \quad \quad \quad \! \qquad \! \qquad j\in S \label{eq:SDPGenEn}
\end{align}

The complexity of problem (\ref{eq:SDPGenSt}--\ref{eq:SDPGenEn}) essentially stems from two fronts, namely (i) computing an expectation, and (ii) finding an optimal policy $\Pi^*(s_{t}, p_{t-1},w_{t-1})$. The complexity of (i) is related to the dimension of the random vector $\tilde{W}_{t-1}$, and step (ii) may be prohibitive because of the dimension of the joint state space $(s_{t},p_{t-1},\tilde{W}_{t-1})$, except for very rare cases, such as \textit{linear systems with quadratic costs} where the optimal policy is a linear function of the state variables (see \cite[pp. 148-149] {bertsekas2005dynamic}).

\section{Approximation through linear programming}\label{GLP}

We will later implement a stochastic dynamic programming scheme to approximate the cost-to-go function. This method assumes linear cost functions, which, however, are usually non-linear. Non-linearity issues may be circumvented via inner \textit{generalized linear programming} (GLP) (see \cite{shapiro1979mathematical}). This is a technique to perform interpolations over a sample of points. In the case  of convexity, GLP converges to the original problem for well conceived sample of points. In case of non-convexity, a convex approximation is performed for the original problem. For each generator $g \in G$, suppose we have exact evaluations of the cost function over the sample of points $\left\{\hat{p}_{gj}|j \in \Lambda_g\right\}$. Assume that $\mathbb{E}_{\tilde{W}_{t}|w_{t-1}}\left[\tilde{F}_{t+1} (s_{t+1},p_t,\tilde{W}_{t})\right]$ is a linear function (we will later substitute the expectation with a single decision variable), the following is a linear approximation to the original problem (\ref{eq:SDPGenSt}--\ref{eq:SDPGenEn}):
\begin{align} \tilde{F}_t(s_{t},p_{t-1},w_{t-1})= & \underset{P_t,\Delta_t,\kappa_t,s_{t+1},\beta_t}{\min{}}\left\{ \sum_{g=1}^{|G|}\sum_{j \in \Lambda_g}  \beta_{gj,t}CP_{g} (\hat{p}_{gj})+\sum^{|S|}_{j=1}CS_{jt}(\Delta_{jt})+ M \sum ^N_{n=1}(\kappa^+_{nt}+\kappa^-_{nt})+ \right. \notag \\
& \quad \qquad \ \ \quad\mathbb{E}_{\tilde{W}_{t}|w_{t-1}}\left[\tilde{F}_{t+1}(s_{t+1}, p_t,\tilde{W}_{t})\right]\label{eq:Linaproxdeb}\bigg\}\\
 & \quad \text{S.t.} \eqref{eq:SDPGenSt}-\eqref{eq:SDPGenEn}\\
& \qquad \quad p_{gt}=\sum_{j\in \Lambda_{g}}\beta_{gj,t}\hat{p}_{gj}, \qquad \quad \qquad g \in G \label{eq:LinaproxProd}\\
& \qquad \quad \sum_{j\in \Lambda_{g}} \beta_{gj,t}=1, \qquad \ \quad \quad \qquad g \in G \label{eq:convexcomb}\\
& \qquad \quad \beta_{gj,t}\geq 0, \qquad \qquad \quad \quad \qquad g \in G, j \in  \Lambda_{g} \label{eq:Linaproxend}
\end{align}

In each period $t$, for each generator $g$, Eq. \eqref{eq:LinaproxProd}, interpolates the production over the sample of generating points $\left\{\hat{p}_{gj}|j \in \Lambda_g\right\}$ using convex coefficients $\beta_{gj,t}$ as computed in Eqs. \eqref{eq:convexcomb} and \eqref{eq:Linaproxend}. Consequently, in the objective, we seek the best interpolation for the cost functions. Though in general, we cannot foresee in which direction GLP errs, however, in the convex case, GLP overestimates the cost functions.

\subsection{Illustration of GLP}
Let us illustrate GLP on the following small quadratic problem:
\begin{align*} 
f(x)= &\underset{x_1,x_2}{\min{}}\left\{\frac{1}{2}x_1^2 + x_2^2 - 2x_1 -6 x_2\right\}\\
& \quad \text{S.t. } x_1 + x_2 \leq 1.75\\
& \qquad \quad - x_1 + 2x_2 \leq 2.2\\
& \qquad \quad 2x_1 + x_2 \leq 4.7
\end{align*}

Solving this problem with Cplex 12.5.1.0, the optimal solution is $x^*_1=0.43333$, $x^*_2 =1.31667$, and the optimal value is $f^*(x)=-6.939167$. Table \ref{tab:GLP} shows the implementation of GLP for different  sample points. GLP quickly converges to the true optimal solution as the grid points is slightly densified. Also observe that the function value is always overestimated.

\begin{table}[!h]
\centering
\caption{Example of GLP implementation}
\label{tab:GLP}
\begin{tabular}{|cc|cc|}\hline
$\hat{x}_1$ & $\hat{f}(x_1)$ & $\hat{x}_2$ & $\hat{f}(x_2)$\\ \hline
3 & -1.5&7&7\\
31&418.5&21&315\\ \hline
$x^*_1$ & $x^*_2$ & $\hat{f}^*(x)$ & \text{True } $f(x)$\\ \hline
-- & -- &  -- & --\\ \hline
$\hat{x}_1$ & $\hat{f}(x_1)$ & $\hat{x}_2$ & $\hat{f}(x_2)$\\ \hline
3 & -1.5&7&7\\
31&418.5&21&315\\ 
0 & 0.0 & 0 &0.0\\ \hline
$x^*_1$ & $x^*_2$ & $\hat{f}^*(x)$ & \text{True } $f(x)$\\ \hline
1.75 & 0.0 &- 0.875 & -1.96875\\ \hline
$\hat{x}_1$ & $\hat{f}(x_1)$ & $\hat{x}_2$ & $\hat{f}(x_2)$\\ \hline
3 & -1.5&7&7\\
31&418.5&21&315\\ 
0 & 0.0 & 0 &0.0\\ 
9&22.5&5&-5\\ \hline
$x^*_1$ & $x^*_2$ & $\hat{f}^*(x)$ & \text{True } $f(x)$\\ \hline
0.43333 & 1.31667 & -1.5333& -6.93917\\ \hline
$\hat{x}_1$ & $\hat{f}(x_1)$ & $\hat{x}_2$ & $\hat{f}(x_2)$\\ \hline
3 & -1.5&7&7\\
31&418.5&21&315\\ 
0 & 0.0 & 0 &0.0\\ 
9&22.5&5&-5\\ \hline
0.5 & -0.875 & 1.5 & -6.75\\\hline
1.75&-1.96875&1.25&-5.9375\\\hline
$x^*_1$ & $x^*_2$ & $\hat{f}^*(x)$ & \text{True } $f(x)$\\ \hline
0.43333 & 1.31667 & -6.5& -6.93917\\ \hline
\end{tabular}
\end{table}
\section{Approximation of the cost-to-go function}\label{sec:ASDP}
Problem (\ref{eq:Linaproxdeb}--\ref{eq:Linaproxend}) cannot be solved for all state values $(s_{t},p_{t-1},w_{t-1})$. A common approach is to discretize $s_{t}, p_{t-1}$, and $\tilde{W}_{t-1}$ spaces into ``partial grids'', and to solve the problem over the Cartesian product of these grids. For instance, suppose we have a network with five wind turbines, five storage devices and five generators. In each period, if we discretize the wind turbine outputs, the production, and the storage into ten level each, without any consideration regarding the computational effort to optimize and to take the expectation, we should solve the problem over a $10^5\times10^5\times10^5=10^{15}$-level grid, which is impossible in practice.  Since its inception, dynamic programming (DP) has been limited to small instance of problems because of the  ``curse of dimensionality'' as coined by its author, Bellman. The computation burden increases exponentially as the number of states increases.

These observations may have significant impact on the practicality of SDP as a solution method for real world problems. This is reinforced by the fact that in certain circumstances, the problem is to be solved periodically within constrained time frame. For instance, in the case of the unit commitment/economic dispatch problem, the operator of the system may need to adjust generating decisions based on signals from the market, or may have to resolve the problem periodically as a result of data changes or updated forecast as he/she gets to receive new observations for the random variables.

The practical limitation of DP paved the way for approximate dynamic programming (ADP) schemes. The purpose is usually to strike a balance between solution time and reasonable performance of the prescribed policy by replacing function $\tilde{F}_t(s_{t},p_{t-1},w_{t-1})$ with some approximation $\hat{F}_t(s_{t},p_{t-1},w_{t-1})$. Though he did not refer to as ADP, Bellman was the first to propose approximations of what he called the \textit{functional equation} (e.g., here Eq. \eqref{eq:SDPGenSt}). Using Lagrange relaxation technique, \cite{bellman1956dynamic}  (see also \cite{bellman1962applied}) discusses successive approximations (SA) of the functional equation  by partitioning the computation of the original sequence of functions into the computation of a sequence of functions of fewer state variables. \cite{korsak1970dynamic} provide a detailed algorithm for the SA technique, and settle conditions for the convergence to true optimal solution. Examples of applications of such technique to power production are available in \cite{zurn1975generator, turgeon1980optimal, yang1989short}. 

As indicated earlier, any DP implementation resorts to discretization of the state space. One straightforward approximation of the cost-to-go function is to select, in each period, a small sample of states instead of a dense grid, and solve problem (\ref{eq:Linaproxdeb}--\ref{eq:Linaproxend}) for each point of the sample, assuming that the minimization as well as the expectation computation can be carried out efficiently. Then the cost-to-go function may be approximated for any out-of-grid point, for instance by interpolation (e.g., linear, multilinear) of the neighboring grid states (see \cite{tejada1993comparison, johnson1993numerical, keane1994solution, philbrick2001improved}). Based on concavity assumptions on the cost-to-go function and assuming that the state space is a hyperrectangle, \cite{zephyr2015controlled} propose a simplicial approximation scheme guided by local estimations of the approximation error. However, the complexity of hypercube decomposition limits the scope of the method.

Several approximation schemes fall within the broader class of \textit{parametric approximation}, whereby the approximation may be denoted by $\hat{F}_t(s_{t},p_{t-1},w_{t-1},\delta)$; $\delta$ is a set of parameters or weights, which usually are to be determined. For instance, in \cite{bellman1957some} such parameters are  Legendre polynomial coefficients. Other polynomial types traditionally used to approximate the cost-to-go function include orthogonal, Chebyshev, spline, Hermite polynomials. For more account on polynomial approximations, see \cite{johnson1993numerical, rust1996numerical, howitt2002using, topaloglu2006dynamic}.

\cite{philbrick2001improved} perform interpolations  for out-of-grid states within hypercubes using weighted sums of the cost-to-go function evaluations as well as derivatives (first and second order) of the function at the vertices. The weights are defined as polynomials in the state variables. This multi-dimensional interpolation approach stems from previous works by \cite{kitanidis1986hermite, kitanidis1987first}, and is an extension of the ``gradient dynamic programming'' scheme by \cite{foufoula1988gradient}.

Combining ideas from various fields such as  \textit{neural networks}, \textit{artificial intelligence}, cognitive sciences, and so on, \textit{reinforcement learning} (or \textit{neuro-dynamic programming}) iteratively constructs approximations to the optimal cost-to-go function, or its expected value \cite{bertsekas1995neuro} through simulations. This methodology is applied by \cite{johri2011self, ahamed2002reinforcement,  vlachogiannis2004reinforcement, momoh2005unit}, to power control decision problems.

 Though there exist several reinforcement learning techniques, the most popular is the \textit{Q-learning} algorithm, which in contrast with DP, computes the cost-to-go function for a set of randomly selected decisions considering only ``visited states''. For further details on contextual applications of this algorithm to power system problems, see \cite{xiong2002electricity, naghibi2006application, lee2007stochastic, tan2009adaptive}.

Reinforcement learning (neuro-dynamic programming) is built around \textit{policy iteration}, and \textit{value iteration}, two widespread DP algorithms \cite{bertsekas1995neuro}. The former algorithm iteratively alternates between \textit{policy evaluation} and \textit{policy improvement}, until no further improvement can be achieved (see \cite{quadrat1988dynamic, qiu1999dynamic, anderson2011adaptive}). The latter algorithm is a later name in the DP literature for successive approximation (see \cite{song2000optimal, anderson2011adaptive}). 

A complete review of ADP methods is beyond the scope of this paper. 

\section{Approximation through stochastic dual dynamic programming} \label{sec:SDDP}
Stochastic dual dynamic programming (SDDP) originated from the work by \cite{pereira1985stochastic, pereira1989optimal,  pereira1991multi}, and may be viewed as a variant of ADP techniques \cite{shapiro2013risk}, where the cost-to-go function is approximated trough sampling and temporal decomposition, alternating between forward and backward steps. 

Unlike other ADP methodologies, in SDDP, the state space is not discretized, but rather sampled. Then the cost-to-go function is approximated in the neighboring of the sampled states by supporting hyperplanes akin to \textit{Benders' cuts} (see \cite{benders1962partitioning}). As a consequence, in principle, SDDP cannot be efficiently implemented for non-convex problems. Several extensions have been proposed to handle non-convexity. For instance, \cite{goor2010optimal, cerisola2012stochastic} replace the non-convex productions functions with piecewise functions, and \textit{McCormick envelopes}, respectively.

Let $V_{t+1}(s_{t+1},p_t,\tilde{W}_{t})=\mathbb{E}_{\tilde{W}_{t}|w_{t-1}}[\tilde{F}_{t+1}(s_{t+1},p_t,\tilde{W}_{t})]$. In each period, SDDP iteratively constructs an ``outer'' approximation to function $V_{t+1}(s_{t+1},p_t,\tilde{W}_{t})$,

\[\hat{V}_{t+1}(s_{t+1},p_t,\tilde{W}_{t}):=\underset{
i}{max}\{H_{t+1}^i(s_{t+1},p_t,\tilde{W}_{t})| i \in I\},\] where $H_{t+1}^i(s_{t+1},p_t,\tilde{W}_{t})$ is a hyperplane constructed using primal and dual information provided by the optimal solution. $\hat{V}_{t+1}(s_{t+1},p_t,\tilde{W}_{t})$ is therefore a lower bound on the true expected cost-to-go $V_t(s_{t},p_{t-1},\tilde{W}_{t-1})$. Figure \ref{fig:SDDPAp} illustrates such an approximation for an hypothetical case.

\begin{figure}[!!h]
\centering
\begin{tikzpicture}[scale =0.7]
\draw[->](0,0)--(0,8); \draw[->](0,0)--(9,0); \draw(10,0) node{$\text{State}$}; \draw(0,8.5) node{$\text{Cost-to-go}$};
\draw  [domain=1:9]  plot (\x, 8/\x^2);  
\draw(4.5,0.)--(0,2.6667); \draw(1,0.5926)--(9,.0006); \draw(0.5625,8)--(2.25,0);
\draw[<-](1.1,8)--(1.5,8);  \draw (2,8) node{$V_{t+1}$};
\draw[<-](1.05,6.125)--(1.45,6.125);  \draw (1.9,6.125) node{$\hat{V}_{t+1}$};
\draw[->](0.5,3)--(0.5,2.5);  \draw (0.7,3.3)node{$H_{t+1}^1$};
\draw[->](1.4,1)--(2,1);  \draw (1,1.3)node{$H_{t+1}^2$};
\draw[->](0.5,0.3)--(1.2,.5); \draw (0.5,0.3)node{$H_{t+1}^3$};
\draw[very thick](0.5625,8)--(1.9685,1.5)--(3.8567,0.3812)--(9,.0006);
\end{tikzpicture}
\caption{Illustration of the cost-to-go function approximation by SDDP}
\label{fig:SDDPAp}
\end{figure}
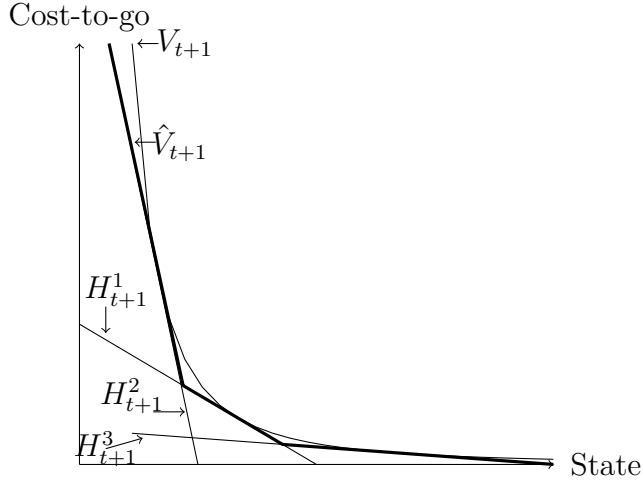

Recall that for any function $f$ of three variables $x, y$, and $z$, the  first-order Taylor approximation, $\hat{f}$, in the neighboring of the point $(x_0,y_0,z_0)$ is given by $\hat{f}(x,y,z)=f(x_0,y_0,z_0)+g_x'(x-x_0)+g_y'(y-y_0)+g_z'(z-z_0)$, which reduces to $\hat{f}(x,y,z)=c+g_x' x + g_y' y+g_z' z$, where $g_x \in \partial_x f(x_0,y_0,z_0)(g_y\in \partial_y f(x_0,y_0,z_0), \text{resp. } g_z \in \partial_z f(x_0,y_0,z_0))$ is a partial subgradient vector of $f$ at $x(y, \text{resp. } z)$, and 
\begin{equation} c=f(x_0,y_0,z_0)-g_x 'x_0-g_y 'y_0-g_z' z_0\label{eq:const}.\end{equation} For any vector $a$, $g_a$ is a column vector, and $g_a'$ is its transpose. We will henceforth drop the transposition operator for simplicity.

 For a given set of ``trial points'' $\Theta_{t+1}=\left\{\left(s_{t+1}^i,p^i_t, w_t^i\right)|1\leq i\leq I\right\}$, we then have 
\[H_{t+1}^i(s_{t+1},p_t, w_t)=\tilde{c}_{t+1}^i+\tilde{g}_{s_{t+1}^i}s_{t+1} +\tilde{g}_{p^i_{t}}p_t + \tilde{g}_{w_t^i} w_t, \ 1\leq i\leq I.\]

$\tilde{c}$ and $\tilde{g}$ are expected values with respect to the random variable $\tilde{W}_{t+1}$.  It is readily verified that
\[\hat{V}_{t+1}(s_{t+1},p_t,w_{t}):=\min{\left\{\rho_{t+1}|\rho_{t+1}\geq \tilde{c}_{t+1}^i+\tilde{g}_{s_{t+1}^i}s_{t+1} +\tilde{g}_{p^i_{t}}p_t + \tilde{g}_{w_t^i} w_t , \ 1\leq i\leq I\right\}}. \]

As a result, an approximation to function $\tilde{F}_t(s_{t},p_{t-1},w_{t-1})$ reads:
\begin{align} \hat{F}_t(s_{t},p_{t-1},w_{t-1})= & \underset{P_t,\Delta_t,\kappa_t,s_{t+1},\rho_{t+1},\beta_t}{\min{}}\left\{\sum_{g=1}^{|G|} h_t(\hat{p}, \Delta_t, \kappa_t)+ \rho_{t+1}\label{eq:aproxdeb}\right\}\\
 & \quad \text{S.t.} \eqref{eq:Linaproxdeb}-\eqref{eq:Linaproxend}\\
& \qquad \quad \rho_{t+1}-\tilde{g}_{s_{t+1}^i}s_{t+1}-\tilde{g}_{p^i_{t}}p_t\geq \tilde{c}_{t+1}^i+ \tilde{g}_{w_t^i} w_t, \qquad \quad 1\leq i\leq I, \label{eq:aproxend}
\end{align}
where $h_t(\hat{p}, \Delta_t, \kappa_t)=\sum_{j \in \Lambda_g}  \beta_{gj,t}CP_{g} (\hat{p}_{g})+\sum^{|S|}_{j=1}CS_{jt}(\Delta_{jt})+ M\sum^N_{n=1}(\kappa_{nt}^+ + \kappa_{nt}^-) $.

\subsection{Building the approximations}

In period $t$, we solve problem (\ref{eq:Linaproxdeb}--\ref{eq:Linaproxend}) for a discrete set of state values $(s_{t}, p_{t-1}, w_{t-1})$ using an approximation, $\hat{V}_{t+1}(s_{t+1}, p_t, w_t)$, of the expected cost-to-go ${V}_{t+1}(s_{t+1}, p_t, w_t)$, defined as a set of $I$ additional constraints to the problem (see problem (\ref{eq:aproxdeb}--\ref{eq:aproxend})). These constraints are built using cut parameters ($\tilde{g}_{s_{t+1}}, \tilde{g}_{p_t}, \tilde{g}_{w_t}$, and $\tilde{c}_{t+1}$) computed in period $t+1$. To account for the uncertainty on the wind turbine outputs, suppose that in period $t+1$, we have considered a finite discrete set of realizations  $\Omega_{t+1}=\left\{w_{t+1}^{j}| 1\leq j\leq J \right\}$. Therefore, in period $t+1$, for each state value $(s_{t+1}, p_t, w_t)$, the problem is solved for each $w_{t+1} \in \Omega_{t+1}$. For each such value of the random process, we have the following storage evolution equation: 
\[s_{t+2}-\Delta^+_{t+1}+\Delta^-_{t+1}=\alpha s_{t+1}.\] Let $\pi_{s,t+1}^j$ be the vector of dual multipliers associated with each such constraint. An expected partial subgradient $\tilde{g}_{s_{t+1}}$ may then be computed as
\[\tilde{g}_{s_{t+1}}=\alpha \sum_{j=1}^J \omega_j \pi_{s,t+1}^{j},\] where $\omega_j$ is the probability of the observation $w^{j}_{t+1}$, and $\sum_{j=1}^J \omega_j=1$.

In period $t+1$,  for each generator $ g\in G$, for each $w_{t+1} \in \Omega_{t+1}$, we also have 
\[ p_{g,t+1} \geq \max\{p_{gt}-\underline{\lambda}_g,\underline{p}_g\},\] 
\[p_{g,t+1} \leq \min \{p_{gt}+\bar{\lambda}_g, \bar{p}_g\}.\]
Let $\pi_{gl}^{j}$ and $\pi_{gu}^{j}$ be the dual prices associated with these constraints, respectively. We know from duality theory (theorem of complementary slackness) that those multipliers are respectively non-null, if and only if the constraints are binding. Consequently, let
\begin{equation*}
g_{gl}^{j}=\left\{ \begin{array}{ll}
\pi_{gl}^{j} & \text{if }p_{g,t+1}^{j*} = p_{gt}-\underline{\lambda}_g,\\
0 & \text{otherwise},
\end{array} \right.
\end{equation*}
where $p_{g,t+1}^{j*}$ is the optimal generation of the generator $g$. Similarly, let
\begin{equation*}
g_{gu}^{j}=\left\{ \begin{array}{ll}
\pi_{gu}^{j} & \text{if }p_{g,t+1}^{j* }= p_{gt}+\bar{\lambda}_g,\\
0 & \text{otherwise}.
\end{array} \right.
\end{equation*}

We then take
\[\tilde{g}_{p_t} = \sum_{j=1}^J \omega_{j} (g_{l}^{j} + g_{u}^{j}).\]

The generation of wind turbines located in the same region are correlated because of similar environmental conditions (wind speed). Consequently, suppose the wind generation is modeled as a lag-$p$ multivariate autoregressive process:
\begin{equation}w_{t+1}= \mu+ \sum^{p-1}_{i=0}\Phi_j (w_{t-i}-\mu)+\tilde{\epsilon}_{t+1}\label{eq:Markov},\end{equation} where $\Phi_j, 0 \leq i \leq p-1$,  is an $M\times M$ matrix of lag- $i+1$ autoregressive coefficients, $\mu$ is the mean vector of the process, $\tilde{\epsilon}_{t+1}$ is a vector of innovations. For a given state value $(s_{t+1}, p_t, w_t)$, for each value $w^{j}_{t+1}$ of the random process, and for each bus where a wind farm $m \in M$ is located, we have the power balance equation
\begin{equation} \sum^{N_m}_{k=1}p_{k,t+1}+\sum_{l\in I_m}e_{l,t+1}-\gamma_{m,t+1}-\kappa_{m,t+1}^+ +\kappa_{m,t+1}^- -\sum_{l\in O_m}e_{l,t+1}=\hat{D}_{m,t+1}-w_{m,t+1}^{j} \label{eq:duabal}.\end{equation} Let $\pi_{wd,t+1}^{j}$ be the vector of dual multipliers associated with these constraints. In addition, suppose the expected cost-to-go function $\hat{V}_{t+2}(s_{t+2}, p_{t+1}, w_{t+1})$ is approximated through the following $K$ cuts: 
\begin{equation}\rho_{t+2}-\tilde{g}_{s^k _{t+2}}s_{t+2}-\tilde{g}_{p_{t+1}^k}p_{t+1}\geq \tilde{c}_{t+2}^k+\tilde{g}_{w_{t+1}^k}w_{t+1}^{j}, 1\leq k\leq K, \label{eq:duaap}\end{equation} (see (\ref{eq:aproxdeb}-\ref{eq:aproxend})). Let $\pi_{wc,t+1}^{jk}$ be the dual multiplier associated with inequality $k$. Using the dual prices associated with Eq. \eqref{eq:duabal} and inequalities \eqref{eq:duaap}, respectively, and Eqs. \eqref{eq:Markov},  we compute a partial subgradient as 
 \[g_{w_t^{j}}=\Phi_1\left(- \pi_{wd,t+1}^{j}+\sum_{k=1}^K \pi_{wc,t+1}^{jk}\tilde{g}_{w_{t+1}^k}\right).\] Above equation follows from the chain rule $\frac{\partial \hat{F}_{t+1}}{\partial w_t}=\frac{\partial w_{t+1}}{\partial w_t}\frac{\partial \hat{F}_{t+1}}{\partial w_{t+1}}$. Thus,  an expected partial subgradient $\tilde{g}_{w_t}$  is given by 
\[\tilde{g}_{w_t}=\sum_{j=1}^J\omega_jg_{w_t^{j}}.\]

Lasly, from Eq. \eqref{eq:const}, we see that
\[\tilde{c}_{t+1} =\sum_{j=1}^J \omega_j \hat{F}_{t+1}(s_{t+1},p_t, w_t) - \tilde{g}_{s_{t+2}} s_{t+2}-\tilde{g}_{p_{t+1}}p_{t+1}-\tilde{g}_{w_{t+1}}w_{t+1},\] where $(s_{t+2}, p_{t+1}, w_{t+1})$ is the state observed in $t+2$ when $(s_{t+1}, p_t, w_t)$ was observed in $t+1$.

Contrary to DP, where in each period we use a look-up  table to keep track of the cost-to-go (we keep in memory the values of the discrete states as well as the corresponding cost-to-go), in SDDP, in each period, we only keep in memory the appropriate Lagrange multipliers and the expected optimal value $\hat{F}_{t}^*$, which are passed to the previous period.

The purpose of the backward recursion phase is to construct the approximations of the cost-to-go function,  available in each period and in each iteration in  the form of supporting hyperplanes to the true function. In each period (for each iteration), suppose we have a set of sampled states $\Theta_{t}=\left\{\left(s_{t}^{i},p_{t-1}^i, w_{t-1}^i\right)| 1 \leq i  \leq I\right\}$. For each state vector, suppose we sample $J$ vectors of wind generation $w_{t}^{ij}, 1\leq j \leq J$, each with probability $\omega^j_i=1/J$. Then for each wind output value, we solve the minimization problem (\ref{eq:aproxdeb}-\ref{eq:aproxend}). From the $J$ minimization problems we retrieve the appropriate multipliers to compute the expected value of the parameters $\tilde{g}_{s_{t}^i}, \tilde{g}_{p_{t-1}^i},\tilde{g}_{w_t^i}$,   and $\tilde{c}_{t}^i$. These parameters are used to construct one supporting hyperplane $H_t^i(s_{t}, p_{t-1},w_{t-1})$ for the minimization problem in period $t-1$. Thus, at the end of the recursion, we shall construct $I$ supporting hyperplanes to the expectation of the cost-to-go function to pass to period $t-1$. 

As initial conditions, we set $\Theta_{T+1}=\emptyset$, $\tilde{g}_{s_{T+1}}=0, \tilde{g}_{p_T}=0, \tilde{g}_{w_T}=0$, and $\tilde{c}_{T+1}=0$. A typical backward iteration is summarized in Table \ref{tab:BWSDDP}.
\begin{table}[!!h]
\centering
\caption{Backward induction procedure}
\label{tab:BWSDDP}
\begin{tabular}{l}
\hline \hline
\text{For } $t=T,T-1,\cdots,1$\\
 \qquad \text{Sample state variables }$(s_{t}^{i}, p^i_{t-1}, w_{t-1}^i),1 \leq i  \leq I.$\\
\qquad \text{For } $i=1,\ldots, I$\\
\qquad \qquad \text{Sample $J$ wind generation vectors } $w_{t}^{ij}, \ \ 1\leq j \leq J.$\\
\qquad \qquad \text{For } $j=1,\ldots,J$\\
\qquad \qquad \qquad \text{Solve the minimization problem  (\ref{eq:aproxdeb}-\ref{eq:aproxend})}.\\
\qquad \qquad \qquad \text{Store the values of the appropriate multipliers and of} $\hat{F}_t^*$.\\
\qquad \qquad \text{End}\\
\qquad \qquad \text{Construct one supporting hyperplane $H_t^i(\cdot)$} for the problem in period $t-1$.\\
\qquad \text{End}\\
\text{End} \\
\hline\\
\end{tabular}
\end{table}

\subsection{Assessing the quality of the approximations}
In each iteration of the algorithm, we construct approximations to the true expected cost-to-go functions through backward recursions. As for any iterative scheme, after each iteration, we need to assess the quality of the approximations, and decide whether or not we can stop the algorithm. In SDDP, such an assessment step is carried out through forward simulations.

For a given initial state $(s_{1}, p_0, w_0)$, suppose we simulate $L$ series of wind generation vectors $w_t^l, 1\leq t\le T, 1\leq l \leq L$, with probability $\omega_l=1/L$ each. Solving problem (\ref{eq:aproxdeb}-\ref{eq:aproxend}) for the first period using the last approximation to the expected cost-to-go function constructed in the  backward step as well as the simulated wind generations, we obtain a lower bound on the true optimal cost (\ref{eq:ProbGenSt}--\ref{eq:ProbGenEn}). Indeed, as we indicated earlier, by construction, in the case of convexity, we have $\hat{F}_2(s_2,p_1,w_1)\leq \tilde{F}_2(s_2,p_1,w_1)$, which implies $\mathbb{E}_{\tilde{W}_2|w_1}\left[\hat{F}_2(s_2,p_1,w_1)\right]\leq \mathbb{E}_{\tilde{W}_2|w_1}\left[\tilde{F}_2(s_2,p_1,w_1)\right]$. It then follows that $\hat{F}_1(s_1,p_0,w_0)\leq \tilde{F}_1(s_1,p_0,w_0)$. Denote this lower bound by $\underline{F}$. 

The decisions obtained from solving the problem in the first period are feasible, but not necessary optimal. By similarly computing the ``immediate'' decisions (production and storage variation) in each period $t, t=2,\ldots,T$, using the cost-to-go approximations constructed in the backward recursion phase, we obtain a series of feasible, but not necessarily optimal management decisions. Thus, the sum of the immediate cost of the decisions so calculated may be used as an estimated upper bound on the true optimal total cost over the $T$ planning span. Denoting this estimated upper bound by $\overline{F}$, we have
\[\overline{F}=\sum_{l=1}^L\sum^T_{t=1} h_t(p_t^l,\Delta_t^l,\kappa_t^l)\omega_l, \] where $h_t(\cdot)$ is the immediate cost in period $t$.

Let $f_l=\sum_{t=1}^Th_t(p_t^l,\Delta_t^l,\kappa_t^l), 1 \leq l\leq L$. The uncertainty around the estimation of the upper bound may be estimated at
\[\sigma_{\bar{F}}=\sqrt{\frac{\sum_{l=1}^l(f_l-\overline{F})^2}{L-1}}.\]

This estimation may be used to construct a $(1-\alpha)\%$ confidence interval around the true value of the upper bound. Usually, if $\underline{F}$ is within the confidence interval, the algorithm is stopped assuming that the approximation is statistically accurate, otherwise, another backward approximation step is performed by adding new supporting hyperplanes (around the state values sampled during the simulation) to tighten the bounds. Drawbacks of such stopping criterion are discussed by \cite{shapiro2011analysis}. \cite{asamov2015regularized, philpott2012dynamic} instead use a fixed number of iterations as termination criterion. Also observe that in each iteration of the forward simulation, we obtain $l$ series of optimal trajectories for the storage levels as well as the generations. These may be used as sampled states for the next backward recursion step. Table \ref{tab:FWSDDP} summarizes the forward procedure.

\begin{table}[!!h]
\centering
\caption{Forward simulation procedure}
\label{tab:FWSDDP}
\begin{tabular}{l}
\hline
\hline
\text{Select an initial state $(s_1, p_0,w_0)$.}\\
\text{Simulate $L$ series of wind generations each of length $T$}.\\
\text{For $l=1,\ldots,L$}\\
\qquad \text{For $t=1,\ldots,T$}\\
\qquad \qquad \text{Solve the minimization problem  (\ref{eq:aproxdeb}-\ref{eq:aproxend})}.\\
\qquad \qquad \text{Store the optimal decisions $s_{t+1}^{l^*}$ and $p_{t}^{l^*}$ and the wind value $w^l_t$}.\\
\qquad \text{End} \\
\text{End} \\
\hline\\
\end{tabular}
\end{table}

\section{Illustration and comparison with classical DP}\label{sec: sddp-dp-illus}
We used the IEEE 9-bus configuration (Figure \ref{fig:ieee-9}) to illustrate and compare SDDP algorithm with classical DP. This network comprises three conventional generators for a total capacity of $820 \ MW$, and nine transmission lines. We used NYISO scaled average hourly load data for January 2016 \cite{nyload}. We located one wind farm and one storage facility at bus 5. Wind data were obtained from the website of the National Renewable Energy Laboratory \cite{wind9}. We utilized a 15\% wind integration scenario and assumed 30\% wind capacity factor. 

Both SDDP and DP algorithms were implemented in Python 2.7.10, and the GLP problems were solved with Cplex 12.5.1.0.  At each iteration of the SDDP backward pass, we simulated 25 vectors of wind outputs to compute the expected values of the cost-to-go functions as well as the expected values of the  partial subgradients ($\tilde{g}_{s_{t}}, \tilde{g}_{p_{t-1}}, \tilde{g}_{w_{t-1}}$), and the expected values of the intercepts ($\tilde{c}_{t}^i$). We also used 25 series of state values in each iteration of the forward pass. As a result, in each iteration of the backward pass, 25 news cuts were added to the GLP problems.

We considered a hypothetical battery with maximum capacity of $250\ MWh$. Maximum hourly charging or discharging  was assumed to be 30\% of the battery capacity, and the minimum level of the stored energy to be 20\% of the battery capacity. For the sake of illustration, we supposed that the battery had 100\% storage efficiency.

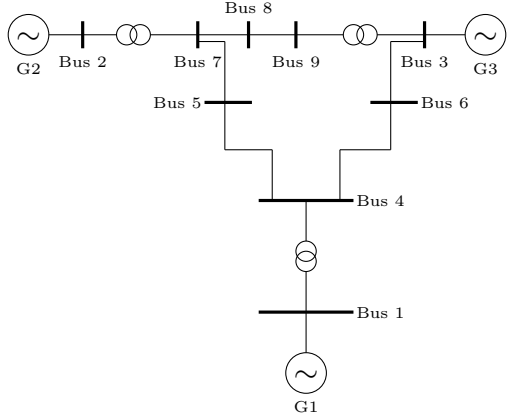
\begin{figure}[h!]
\centering
\begin{tikzpicture}[scale=0.9]
\draw (0,0) circle(0.3); \draw(0.,-0.05) node{$\sim$};  \draw(0,-0.5) node{$\text{\tiny{G2}}$}; \draw (0.3,0.)--(0.8,0); \draw[very thick] (0.8,-0.2)--(0.8,0.2); \draw(0.8,-0.4) node{$\tiny{\text{Bus 2}}$}; \draw (0.8,0.)--(1.3,0); \draw (1.45,0) circle (0.15); \draw (1.65,0) circle (0.15); \draw (1.8,0.)--(2.5,0); \draw[very thick] (2.5,-0.2)--(2.5,0.2); \draw(2.5,-0.4) node{$\tiny{\text{Bus 7}}$}; \draw (2.5,0.)--(3.25,0); \draw[very thick] (3.25,-0.2)--(3.25,0.2); \draw(3.25,0.4) node{$\tiny{\text{Bus 8}}$}; \draw (3.25,0.)--(3.95,0); \draw[very thick] (3.95,-0.2)--(3.95,0.2); \draw(3.95,-0.4) node{$\tiny{\text{Bus 9}}$}; \draw (3.95,0.)--(4.65,0); \draw (4.8,0) circle (0.15); \draw (5.0,0) circle (0.15); \draw (5.15,0.)--(5.85,0);
\draw[very thick] (5.85,-0.2)--(5.85,0.2); \draw(5.85,-0.4) node{$\tiny{\text{Bus 3}}$}; \draw (5.85,0.)--(6.45,0); \draw (6.75,0) circle(0.3); \draw(6.75,-0.05) node{$\sim$}; \draw(6.75,-0.5) node{$\text{\tiny{G3}}$};

\draw (2.5,-0.1)--(2.9,-0.1);  \draw (2.9,-0.1)--(2.9,-1);  \draw[very thick] (2.6,-1)--(3.3,-1); \draw(2.2,-1) node{$\tiny{\text{Bus 5}}$}; \draw (5.85,-0.1)--(5.35,-0.1); \draw (5.35,-0.1)--(5.35,-1); \draw[very thick] (5.05,-1)--(5.75,-1); \draw(6.15,-1) node{$\tiny{\text{Bus 6}}$};

\draw (5.35,-1)--(5.35,-1.7); \draw (2.9,-1)--(2.9,-1.7); \draw (2.9,-1.7)--(3.6,-1.7); \draw (5.35,-1.7)--(4.6,-1.7);
\draw (3.6,-1.7)--(3.6,-2.45); \draw (4.6,-1.7)--(4.6,-2.45); \draw[very thick] (3.4,-2.45)--(4.8,-2.45); \draw(5.2,-2.45) node{$\tiny{\text{Bus 4}}$};
\draw (4.1,-2.45)--(4.1,-3.05); \draw (4.1,-3.2) circle (0.15); \draw (4.1,-3.35) circle(0.15); \draw (4.1,-3.5)--(4.1,-4.1);
\draw[very thick] (3.4,-4.1)--(4.8,-4.1); \draw(5.2,-4.1) node{$\tiny{\text{Bus 1}}$}; 
\draw (4.1,-4.1)--(4.1,-4.7); \draw (4.1,-5) circle(0.3); \draw(4.1,-5.05) node{$\sim$}; \draw(4.1,-5.5) node{$\text{\tiny{G1}}$};
\end{tikzpicture}
\caption{9-bus system configuration}
\label{fig:ieee-9}
\end{figure}

Figure \ref{fig:ssdp-9-cs} depicts the trajectories of the storage facility over a 24 hour span for two cases. In case (a) the stored energy is at the minimum allowable level  in the beginning of the first hour; in case (b) the level of charge is at 60\%  in the beginning of the first hour. We also assumed zero-cost for varying the level of energy in this example. The net load (demand - wind output) is also depicted in each case.

In the first case, the battery is charged progressively in the beginning as the net load is low, then is steady as of the first peak hour, and up to the second peak hour. Then the battery is discharged progressively to its minimum level, since the model considered no ``terminal value'' for the storage. In the second situation, the level of charge is steady until the first peak hour, then the battery is depleted until the end of the horizon.
\begin{figure}[h!]
\centering
\scalebox{0.7}[0.7]{\includegraphics[clip, viewport=0.5in 3.3in 8in 7.6in]{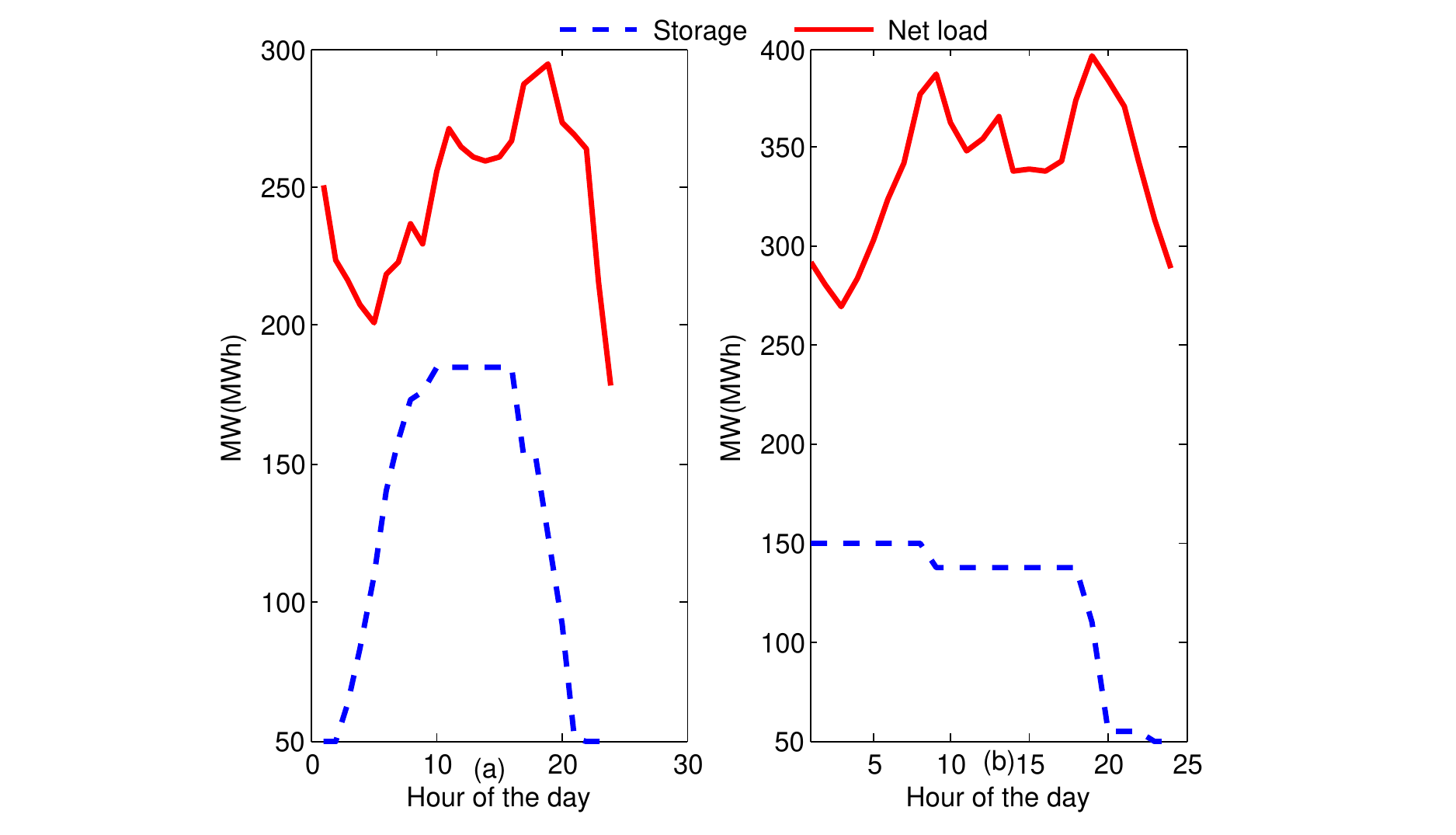}}
\caption{Example of storage trajectory when charging (discharging) cost is low}
\label{fig:ssdp-9-cs}
\end{figure}

Figure \ref{fig:ssdp-9-es} illustrates a situation where varying the level of charge of the battery is very expensive as compared to the cost of operating the conventional units. Since utilizing the battery is very costly, the stored energy is used only in hours where the conventional generators cannot meet the demand.

\begin{figure}[h!]
\centering
\scalebox{0.7}[0.7]{\includegraphics[clip, viewport=0.5in 3.3in 8in 7.6in]{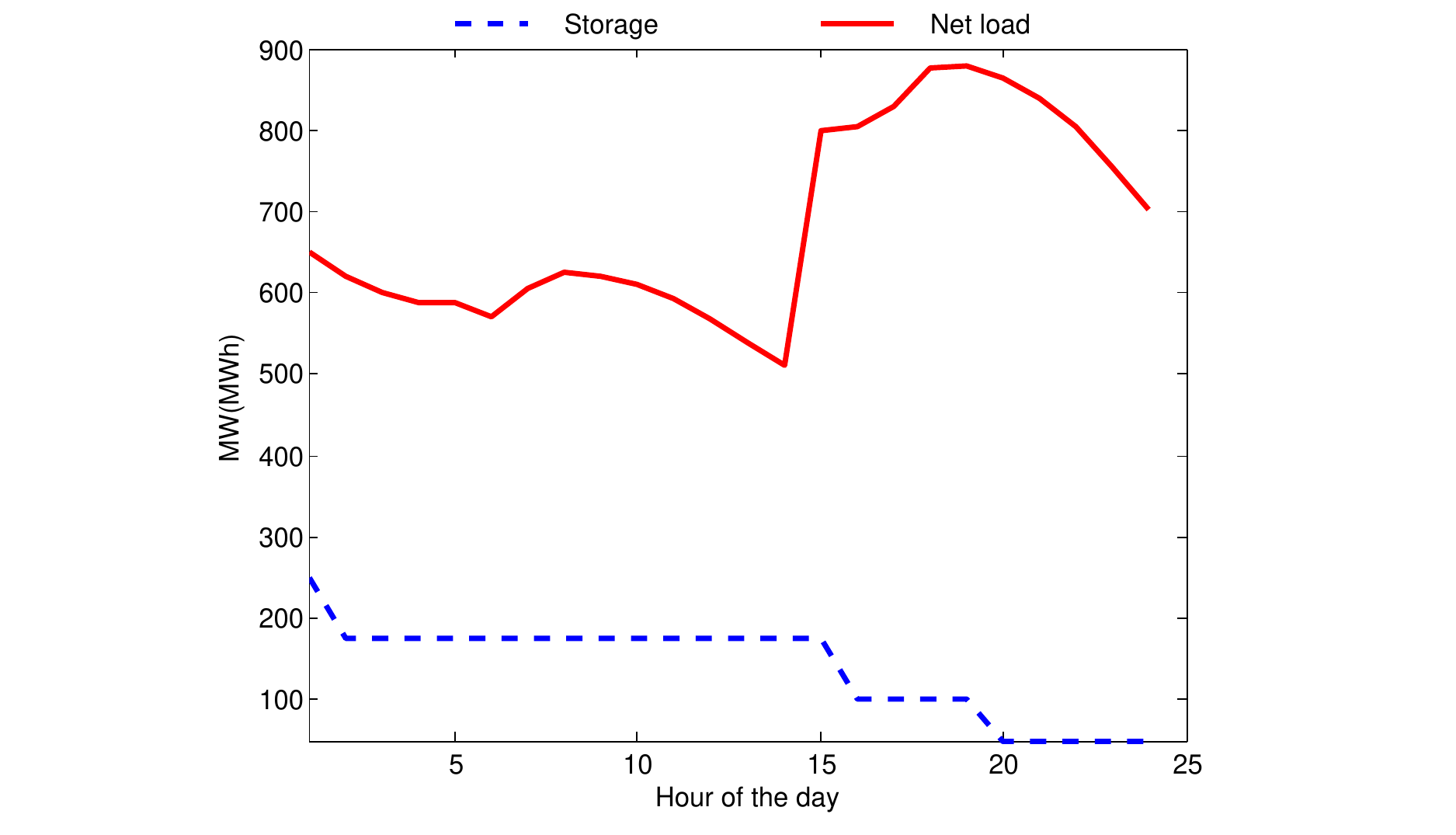}}
\caption{Example of storage trajectory when charging (discharging) cost is high}
\label{fig:ssdp-9-es}
\end{figure}

How would SDDP compare with classical DP? In general, this would difficult to answer due to the  prohibitiveness of the computational burden of DP. We compared the performance of both algorithms on the 9-bus network example based on two criteria, namely (i) the CPU time, and (ii) the solution performance (total cost). We first approximated the value functions with both SDDP and DP. For DP, in each hour, the storage, and the generators outputs were discretized  each into six levels, and the wind farm output into seven levels. This was found to be the finest grids to solve the problems in reasonable time. Since we assumed the wind process to be Markovian, in each hour, we needed to consider all the transitions from the previous wind values (e.g. $7\times 7$ transitions). As a result, in each hour, the DP problem was solved over a grid of $6\times 6^3\times 7^2$ levels (the network comprising three conventional generators), which resulted in 63 504 evaluations of the cost-to-go function in each hour. Therefore, 1 524 096 evaluations were performed over the 24-hour horizon.

Following the construction of the cost-to-go functions with both methods, we simulated the operation of the network for two simulation runs, each comprising one hundred 24-hour wind scenarios.  With  DP, in both run, the cost-to-go functions were constructed over the same state space. With SDDP, in the first run, we performed only four iterations to approximate the value functions. Since in each hour we approximated the expected cost-to-go using 25 supporting hyperplanes (state values) and used 20 wind values to compute the expectations, with SDDP, the total number of function evaluations then were 500, for a total of 12 000 function evaluations over the 24 hour horizon. Therefore, overall 48 0000 function evaluations were performed in the first simulation. Those numbers do not include the SDDP forward step - to assess the quality of the approximation-, since this is very fast. In the second simulation, we tried to improve the quality of the approximation by carrying out ten iterations of SDDP. Consequently, 120 000 function evaluations were performed in that run.

Table \ref{tab:cpu-sddp-dp} presents the CPU time for each method and each run.  In the first case, the computation time for SDDP was about 4\% of that of DP. In the second case, that proportion was about 15\% as better approximations were sought.     

\begin{table}[h!!!]
\centering
\caption{Comparison of the CPU time in seconds: SDDP and DP}
\label{tab:cpu-sddp-dp}
\begin{tabular}{|c|c|c|}\hline
\text{Method} & \text{Run 1}& \text{Run 2} \\ \hline
\text{SDDP} & 499.35&1 876.22\\ \hline
\text{DP} & 13 038.16&12 726.88 \\ \hline
\end{tabular}
\end{table}

Table \ref{tab:cost-sddp-dp} reports descriptive statistics on the performance (total cost) of each method for each run. In the first case, on average (over the 100 scenarios), using DP, the total cost was only improved by 0.71\% as compared to SDDP. In the second case, on average,the total cost difference decreased to 0.21\%. This suggests that, overall, SDDP allowed for a fair trade-off between solution time and accuracy.

\begin{table}[h!!!]
\centering
\caption{Comparison of solution cost: SDDP and DP}
\label{tab:cost-sddp-dp} 
\begin{tabular}{|c|c|c|c|c|}\hline
\multirow{2}{*}{Method} & \multicolumn{4}{c|}{Run 1}\\ \cline{2-5}
                            &\text{Min}& \text{Max}& \text{Mean} & \text{Standard deviation}\\ \hline
\text{SDDP} & 78 622.29 &162 215.38 & 126 872.83&21 812.53\\ \hline
\text{DP} &75 392.76& 161 872.07 & 125 968.26&22 922.10\\ \hline
\multirow{2}{*}{Method} & \multicolumn{4}{c|}{Run 2}\\ \cline{2-5}
                            &\text{Min}& \text{Max}& \text{Mean} & \text{Standard deviation}\\ \hline
\text{SDDP} & 88 691.21& 164 333.09 & 134 267.42&19 210.77\\ \hline
\text{DP} &85 882.05&164 333.09 & 133 981.13&19 657.52\\ \hline
\end{tabular}
\end{table}

Figure \ref{fig:ssdp-dp-trajec} depicts the average storage trajectory for both method and both simulation runs. In general, with DP, the level of charge is higher, but the overall pattern (time of charging and discharging) is similar. 
\begin{figure}[h!]
\centering
\scalebox{0.7}[0.7]{\includegraphics[clip, viewport=0.5in 3.3in 8in 7.6in]{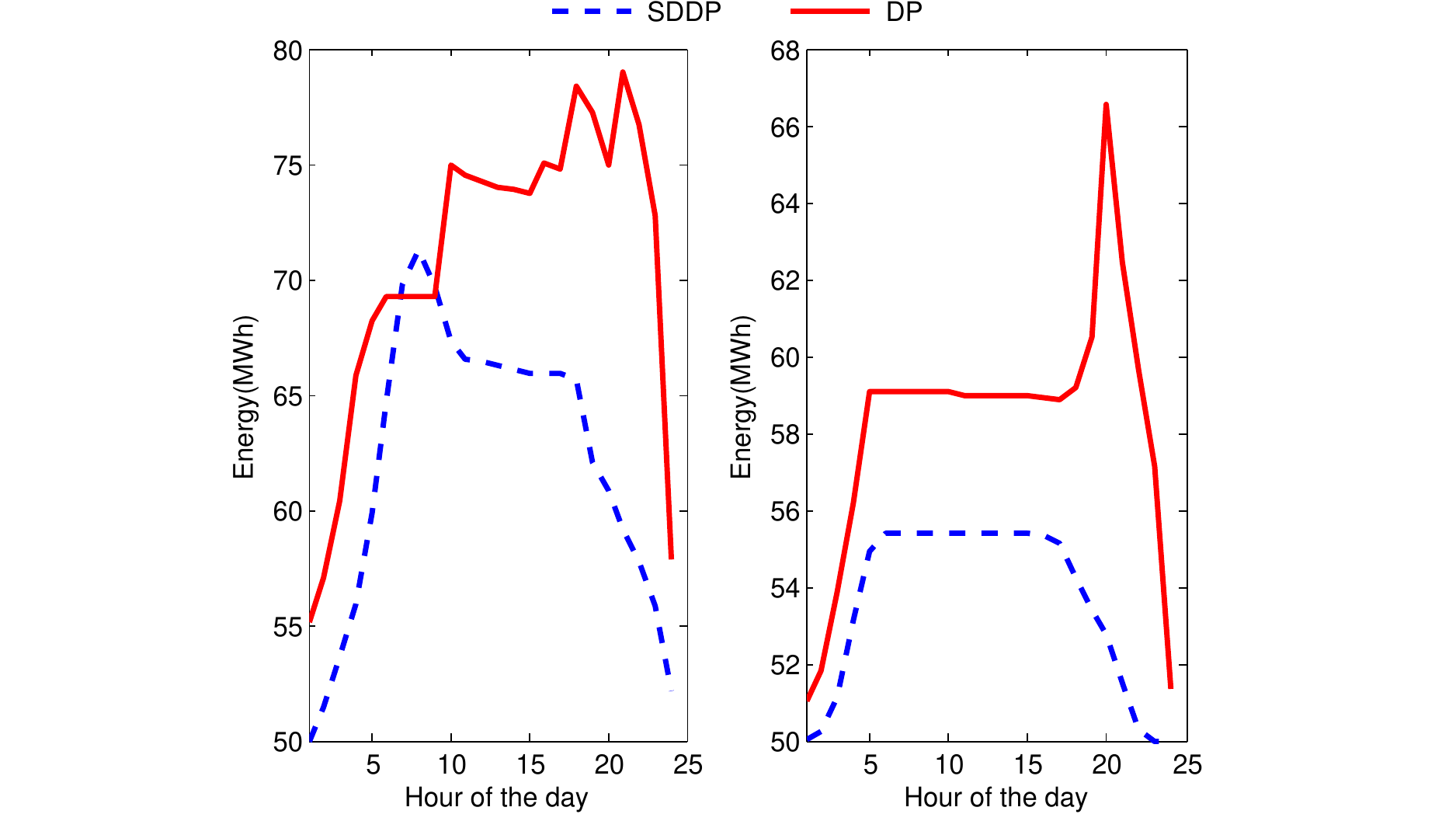}}
\caption{Mean storage trajectory obtained from SDDP and DP}
\label{fig:ssdp-dp-trajec}
\end{figure}

\section{Additional numerical tests}\label{sec:more-num}
The ultimate goal of this research was to analyze the scalability of SDDP to larger networks. We tested the algorithm on different IEEE networks, which characteristics are presented in Table \ref{tab:systemcha}. Each network was tested with different number of storage facilities and wind farms. Wind and load data were obtained from the sources as in Section \ref{sec: sddp-dp-illus}.

\begin{table}[h!!!]
\centering
\caption{Test networks' characteristics}
\label{tab:systemcha} 
\begin{tabular}{|c|c|c|c|}\hline
\text{Case} & \# \text{of buses} & \# \text{of generators} & \# \text{of trans. lines} \\ \hline
1 & 30 & 6 & 41 \\\hline
2 & 57 & 7 & 80  \\\hline
3 & 89 & 12 & 210 \\\hline
4 & 118 & 54 & 186 \\\hline
5 & 300 & 69 & 411\\\hline
\end{tabular}
\end{table}

Figure \ref{fig:trajec-118b-5s} shows an example of the mean trajectory over 100 simulations for a network composed of 118 buses, 5 storage facilities and one aggregated wind farm, as well as the net load. Each battery has specific storage, charging, and discharging efficiency, respectively. The three batteries that either have the highest storage or charging/discharging efficiency are used to contribute to meet the load, whereas the energy level of the other two batteries (with either the lowest storage or charging/discharging efficiency) is steady over the 24-hour span.
\begin{figure}[h!]
\centering
\scalebox{0.7}[0.7]{\includegraphics[clip, viewport=0.5in 3.3in 8in 7.6in]{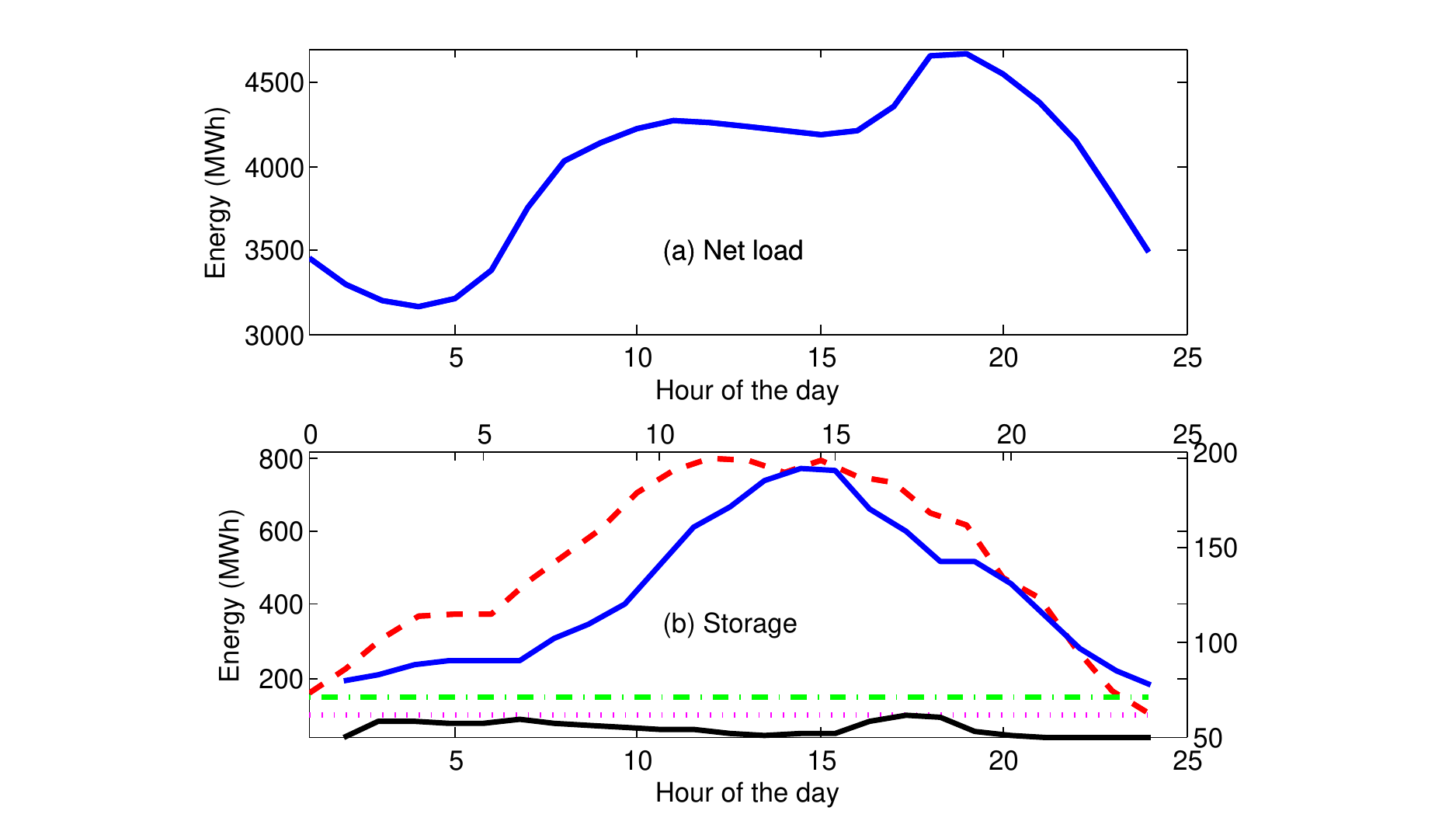}}
\caption{Mean storage trajectory for the 118-bus network with five storage units and one wind farm}
\label{fig:trajec-118b-5s}
\end{figure}

Table \ref{tab:cpu-sddp} reports the computation time for each network, and equipped with different number of storage units (|S|) and wind farms (|M|). We varied the number of storage facilities from one to twenty, and the number of wind farms from one to ten. The computation time varied between: (i) 1 230 seconds and 1 583 seconds in the case of the 30- bus network, (ii) 1 388 seconds and 1 598 seconds for the 57-bus network, (iii)  1 570 and 1 737 seconds for the 57-bus network, (iv) 2 180 seconds and 2 454 seconds in the case of the 118-bus network, and, (v) 4 159.16 seconds and 5 036.37 seconds for the largest network (300 buses).
\begin{table}[h!!!]
\centering
\caption{Computation time in seconds}
\label{tab:cpu-sddp}
\begin{tabular}{|cccc|cccc|} \hline
\# \text{buses} & |S| & |M| & \text{CPU time}&\# \text{buses} & |S| & |M| & \text{CPU time}\\\hline
30 &1 & 1& 1 229.80&118 & 1&1&2 399.35\\\
30 & 5 & 1&1 582.67&118& 5 & 1& 2 444.99 \\
30& 5 & 5 &  1 323.88&118& 5 & 5 & 2 453.89\\
57&1 & 1&1 388.09&118&10&5&2 179.62 \\
57 &5& 1&1 454.47 &118&20&10&2 248.39\\
57&5 & 5 &  1 396.26&300&1&1&4 159.16\\
57&10	&5&1 597.71&300&5&1&4 234.72\\
89 & 1& 1&1 570.67&300&5&5&4 570.01\\
89 & 5 & 1&1 709.68&300&10&5&4 617.65 \\
89& 5 & 5 & 1 575.09&300&20&10&5 036.37\\ 
89&10	&5& 1 737.32&&&&\\ \hline
\end{tabular}
\end{table}

\section{Conclusions}

This paper analyzed the operation of power networks comprising both conventional and wind generators in the presence of storage. We formulated the problem under the framework of dynamic programming, and used stochastic dual dynamic programming and generalized linear programming to approximate the problem. Such approximation schemes allowed to handle large dimension state space and to solve the problem in reasonable time as compared to classical dynamic programming.

\bibliography{Wind_ManagementWithStorage.bbl}
\end{document}